\input amstex
\documentstyle{amsppt}

\document
\magnification=1200
\NoBlackBoxes
\nologo

\vsize18cm

\centerline{\bf Bihamiltonian elliptic structures}
\bigskip
\centerline{\bf Alexander Odesskii}
\bigskip
\centerline{\bf Introduction}
\medskip

Bihamiltonian structures play an important role in the theory of dynamic systems. In this approach one starts 
with two Poisson brackets $\{ , \}_1$ and $\{ , \}_2$ on some manifold, such that any linear combination
 $\{ , \}_{\lambda_1,\lambda_2}=\lambda_1\{ , \}_1+\lambda_2\{ , \}_2$ is also a Poisson bracket. Basing on
 these brackets one constructs a Hamiltonian system. The construction of the dynamical system basing on these
 brackets is called Lenard scheme [1,2,3,4]. It provides a family of functions in involution. Namely, let 
 $C_{\lambda_1,\lambda_2}$ and $C^{\prime}_{\lambda_1,\lambda_2}$ be central elements for the Poisson structure 
 $\{ , \}_{\lambda_1,\lambda_2}$, then one has $\{C_{\lambda_1,\lambda_2},C^{\prime}_{\mu_1,\mu_2}\}_i=0$ for 
 $i=1,2$.
 
 In this paper we construct three Poisson structures $\{ , \}_i$ $(i=1,2,3)$ on $\Bbb C^n$ such that any $\{ , \}_i$ 
 is quadratic and any linear combination $\{ , \}_{\lambda_1,\lambda_2,\lambda_3}=\sum_i\lambda_i\{ , \}_i$ 
 is also a Poisson bracket. We also study symplectic leaves of the Poisson structure
 $\{ , \}_{\lambda_1,\lambda_2,\lambda_3}$ and construct central elements.
 
 Let $\Cal E=\Bbb C\diagup\Gamma$ be an elliptic curve and $\eta\in \Cal E$. The algebra $Q_n(\Cal E,\eta)$ 
 is generated by $n$ elements $\{x_i; i\in \Bbb Z\diagup n\Bbb Z\}$ with the following defining relations:
 
 $$\botshave{\sum_{r\in \Bbb Z\diagup n\Bbb Z}}\frac{1}{\theta_{j-i-r}(-\eta)\theta_r(\eta)}x_{j-r}x_{i+r}=0,$$ 
 for all $i,j\in \Bbb Z\diagup n\Bbb Z$ such that $i\ne j$. Here $\{\theta_i(z), i\in \Bbb Z\diagup n\Bbb Z\}$ 
 are $\theta$-functions of order $n$ with respect to the lattice $\Gamma\subset \Bbb C$. It is known [5,6] 
 that for generic $\eta$ the algebra $Q_n(\Cal E,\eta)$ has the same size of graded components as the 
 polynomial algebra in $n$ variables. Hence, for any fixed elliptic curve $\Cal E$ one has the flat deformation
 of the polynomial algebra. Let $q_n(\Cal E)$ be the corresponding Poisson structure on $\Bbb C^n$. It turns out 
 that the Poisson structure $\{ , \}_{\lambda_1,\lambda_2,\lambda_3}$ for generic point $(\lambda_1,\lambda_2,\lambda_3)
 \in \Bbb C^3$ is isomorphic to $q_n(\Cal E)$ for certain $\Cal E$. It follows that the algebras
 $Q_n(\Cal E,\eta)$ are quantizations of the Poisson structures $\{ , \}_{\lambda_1,\lambda_2,\lambda_3}$ for generic 
  $\lambda_1,\lambda_2,\lambda_3$ and one may expect that there is a lot of commuting elements in the algebras
  $Q_n(\Cal E,\eta)$, which are quantization of commuting elements in $\{ , \}_{\lambda_1,\lambda_2,\lambda_3}$ 
  provided by Lenard scheme.
  
  Now we describe the contents of the paper. In {\bf \S1} we construct three compatible quadratic Poisson structures. 
  In {\bf \S2} we give another description of the Poisson structures in terms of elliptic functions. This construction 
  is called functional realization. In {\bf \S3} we study symplectic leaves of our Poisson structure and construct
  Casimir elements. In {\bf Appendix} we collect some standard facts from the theory of elliptic functions [7].
  
  \newpage
  \centerline{\bf \S1. Three Poisson structures}
  \medskip
  
  For any $n\in\Bbb C$ we define three quadratic Poisson structures in the polynomial algebra with infinite number 
  of generators $\{e_{\alpha}; \alpha\in\Bbb Z\}$. But it will be clear that for $n\in\Bbb N$ the polynomial algebra
  generated by $e_0$ and $\{e_{\alpha}; 2\leqslant\alpha\leqslant n\}$ is Poisson subalgebra for all our Poisson structures.
  We will use the notation: $S_k(e_{\alpha},e_{\beta})=\botshave{\topshave{\sum_{r=0}}^{\infty}}e_{\alpha+kr}e_{\beta-kr}$. 
  It is clear 
  that if $\alpha+\beta=\alpha^{\prime}+\beta^{\prime}$ and $\alpha\equiv\alpha^{\prime}\mod k$, then the linear combination 
  $S_k(e_{\alpha},e_{\beta})-S_k(e_{\alpha^{\prime}},e_{\beta^{\prime}})$ contains only finite number of monomials. 
  Define $\{,\}_i$, $i=1,2,3$, by the following formulas:
  $$\{e_{\alpha},e_{\beta}\}_1=\frac{n}{2}(S_1(e_{\alpha+1},e_{\beta})-S_1(e_{\beta+1},e_{\alpha}))+(\alpha-n)e_{\alpha+1}
  e_{\beta}-(\beta-n)e_{\alpha}e_{\beta+1}$$
  $$\{e_{2\alpha},e_{2\beta}\}_2=0$$
  $$\{e_{2\alpha},e_{2\beta+3}\}_2=\frac{n}{8}(S_2(e_{2\beta+2},e_{2\alpha-2})-S_2(e_{2\alpha},e_{2\beta}))+\frac{1}{4}
  (2\beta+1)e_{2\alpha}e_{2\beta} $$
  $$\{e_{2\alpha+3},e_{2\beta+3}\}_2=\frac{n}{4}(S_2(e_{2\beta+2},e_{2\alpha+1})-S_2(e_{2\alpha+2},e_{2\beta+1}))-$$
  $$-\frac{1}{4}
  (2\alpha+1)e_{2\alpha}e_{2\beta+3} +\frac{1}{4}(2\beta+1)e_{2\alpha+3}e_{2\beta}\eqno(1)$$
  $$\{e_{2\alpha},e_{2\beta}\}_3=0$$
  $$\{e_{2\alpha},e_{2\beta+3}\}_3=\frac{n}{8}(S_2(e_{2\beta},e_{2\alpha-2})-S_2(e_{2\alpha},e_{2\beta-2}))+\frac{1}{2}
  \beta e_{2\alpha}e_{2\beta-2} $$
  $$\{e_{2\alpha+3},e_{2\beta+3}\}_3=\frac{n}{4}(S_2(e_{2\beta},e_{2\alpha+1})-S_2(e_{2\alpha},e_{2\beta+1}))-\frac{1}{2}
  \alpha e_{2\alpha-2}e_{2\beta+3} +\frac{1}{2}\beta e_{2\alpha+3}e_{2\beta-2}$$
  
  {\bf Proposition 1.} {\it For any $n\in \Bbb N$ these formulas (1) define Poisson structures in the polynomial 
  algebra generated by $e_0$ and $\{e_{\alpha}; 2\leqslant\alpha\leqslant n\}$. Moreover, any linear combination of 
  $\{,\}_i$, $i=1,2,3$ is also a Poisson bracket.}
  
  {\bf Proof} of this proposition is straightforward. It also follows from the functional construction in the next
  paragraph, which explains the elliptic nature of this Poisson structure.
  
  \newpage
  \centerline{\bf \S2. Functional realization}
  \medskip
  
  Let us fix an integral lattice $\Gamma\subset\Bbb C$. We will use the standard notations from the theory of elliptic
  functions, like $\wp(z), \zeta(z), g_2, g_3$ (see Appendix).
  
  Let $\Cal F$ be the space of elliptic functions in one variable with respect to $\Gamma$ and holomorphic outside $\Gamma$.
  For $n\in\Bbb N$ let $\Cal F_n\subset\Cal F$ be the subspace of functions with poles of order $\leqslant n$ on $\Gamma$.
  It is clear that $\dim\Cal F_n=n$. It is known that the functions $\{e_{2\alpha}(z)=\wp(z)^{\alpha}, e_{2\alpha+3}(z)=
  -\frac{1}{2}\wp(z)^{\alpha}\wp^{\prime}(z); \alpha\in\Bbb Z_{\geqslant 0}\}$ form a basis of the linear space $\Cal F$. 
  It is clear that for any $\alpha\in\Bbb N$ the function $e_{\alpha}(z)$ has a pole of order $\alpha$ in $\Gamma$ with 
  residue 1. The functions $e_0$ and $e_{\alpha}, 2\leqslant\alpha\leqslant n$ form a basis of the space $\Cal F_n$. 
  It is clear that the symmetric power $S^m\Cal F$ (resp. $S^m\Cal F_n$) is isomorphic to the space of symmetric 
  elliptic functions in $m$ variables $f(z_1,...,z_m)$ holomorphic outside of the divisors $z_p\in\Gamma, 1\leqslant p
  \leqslant m$ (resp. with poles of order $\leqslant n$ on these divisors).
  
  We construct a bilinear operator $\{,\}: \Lambda^2\Cal F\to S^2\Cal F$ as follows: for $f,g\in\Cal F$ we set
  
  $$\{f,g\}(x,y)=n(\zeta(x-y)-\zeta(x)+\zeta(y))(f(x)g(y)-f(y)g(x))-$$
  
  $$-f^{\prime}(x)g(y)-f^{\prime}(y)g(x)+f(x)g^{\prime}(y)+f(y)g^{\prime}(x) \eqno(2)$$
  
  {\bf Proposition 2.} {\it The formula (2) defines a Poisson structure on the polynomial algebra $S^*\Cal F$.
  If $n\in\Bbb N$, then $\{\Cal F_n,\Cal F_n\}\subset S^2\Cal F_n$ and we have the Poisson structure on the polynomial 
  algebra $S^*\Cal F_n$. Moreover, in the basis $\{e_{\alpha}\}$ we have $\{e_{\alpha},e_{\beta}\}=
  \{e_{\alpha},e_{\beta}\}_1+g_2\{e_{\alpha},e_{\beta}\}_2+g_3\{e_{\alpha},e_{\beta}\}_3$, where $\{,\}_i$, $i=1,2,3$ 
  are defined by (1). }
  
  {\bf Proof} is a simple calculation with Weierstrass functions using identities (5) from Appendix. It is clear that 
  if $f$ and $g$ are elliptic functions, then l.h.s. of (2) is also an elliptic function in two variables.
  Moreover, if $f$ and $g$ are holomorphic outside $\Gamma$ with poles of order $\leqslant n$ on $\Gamma$, then $\{f,g\}$
has the same property. Verification of the Jacobi identity is straightforward. For calculation of this Poisson
  brackets in the basis $e_{\alpha}$ one needs only the identities (5) from Appendix.
  
  \newpage
  \centerline{\bf \S3. Symplectic leaves and Casimir elements}
  \medskip
  
  For $p\in\Bbb N$ we denote by $b_{p,n}$ the Poisson algebra which is spanned by the elements
  $\{f(u_1,...,u_p)\psi_1^{\alpha_1}...\psi_p^{\alpha_p}; \alpha_1,...,\alpha_p\in \Bbb Z_{\geqslant 0}\}$ as a linear space, where
  $u_1,...u_p,\psi_1,...,\psi_p$ are independent variables and $f(u_1,...u_p)$ are elliptic functions in variables
  $u_1,...,u_p$ with respect to the lattice $\Gamma$ holomorphic outside the divisors $u_j\in\Gamma$ and 
  $u_j-u_k\in\Gamma$. Poisson bracket on $b_{p,n}$ is defined as follows:
  
  $$ \{u_{\alpha},u_{\beta}\}=0, \{u_{\alpha},\psi_{\beta}\}=\psi_{\beta}, \{u_{\alpha},\psi_{\alpha}\}=
  -\frac{n-2}{2}\psi_{\alpha},$$
  $$ \{\psi_{\alpha},\psi_{\beta}\}=n(\zeta(u_{\alpha}-u_{\beta})-\zeta(u_{\alpha})+\zeta(u_{\beta}))
  \psi_{\alpha}\psi_{\beta}, \eqno(3)$$
  where $\alpha\ne\beta$.
  
  Let us define a linear map $x_p:\Cal F\to b_{p,n}$ by the formula:
  
  $$x_p(f)=\botshave{\sum_{1\leqslant\alpha\leqslant p}}f(u_{\alpha})\psi_{\alpha} \eqno(4)$$
  There is a unique extension of this map to the homomorphism of commutative algebras $S^*\Cal F\to b_{p,n}$ which
  we also denote by $x_p$.
  
  {\bf Proposition 3.} {\it The map $x_p: S^*\Cal F\to b_{p,n}$ is a homomorphism of the Poisson algebras.}
  
  {\bf Proof.} It is easy to check that $x_p(\{f,g\})=\botshave{\sum_{1\leqslant\alpha,\beta\leqslant p}}\{f,g\}(u_{\alpha},
  u_{\beta})\psi_{\alpha}\psi_{\beta}$ for any $f,g\in\Cal F$. This implies the proposition.
  
  If $2p<n$, then the formulas (3) define nondegenerate Poisson structure on the open set $\{\psi_{\alpha}\ne 0;u_{\alpha},
  u_{\alpha}-u_{\beta}\notin\Gamma\}$. In this case the formula (4) defines symplectic leaves of the Poisson algebra 
  $S^*\Cal F_n$. Central elements of the Poisson algebra  $S^*\Cal F_n$ belong to $\ker x_p$ for $2p<n$, because the 
  algebra $b_{p,n}$ is nondegenerated in this case. One can check that for $p=\frac{n}{2}-1$ for even $n$ (resp. 
  $p=\frac{n-1}{2}$ for odd $n$) the ideal $\ker x_p$ is generated by two elements of degree $\frac{n}{2}$ (resp. 
  by one element of degree $n$). Center of the Poisson algebra $S^*\Cal F_n$ is a polynomial algebra generated
  by these elements. In fact, our construction of homomorphism $x_p$ implies that $\ker x_p$ on the space $S^{p+1}\Cal F_n$
  consists of such elements $f(z_1,...z_{p+1})\in S^{p+1}\Cal F_n$, which are equal to zero on the divisors 
  $z_j-z_k\in\Gamma$, $1\leqslant j<k\leqslant p+1$. This allows us to construct central elements from 
  $\ker x_{\frac{n}{2}-1}$ of degree $\frac{n}{2}$ for even $n$ explicitly.
  
  {\bf Examples:}
  
  {\bf 1. $n=2$.} The Poisson algebra $S^*\Cal F_2$ is commutative, one has $\{e_0,e_2\}=0$ for $n=2$.
  
  {\bf 2. $n=4$.} Let $$C^{(2)}_0=\left| \matrix
                e_0 & e_2\\
                e_2 & e_4
                 \endmatrix \right|, C^{(2)}_1=\left| \matrix
                e_2 & e_3\\
                e_3 &  e_4-\frac{1}{4}g_2e_0
                 \endmatrix \right|+\frac{1}{4}g_3e_0^2$$

  These elements are central in the algebra $S^*\Cal F_4$.
  
  {\bf 3. $n=6$.} Let $$C^{(3)}_0= \left| \matrix
                e_0 & e_2 & e_3 \\
                e_2 & e_4 & e_5 \\
                e_3 & e_5 & e_6-\frac{1}{4}g_2e_2-\frac{1}{4}g_3e_0
                 \endmatrix \right|,$$
                 $$C^{(3)}_1= \left| \matrix
                e_2 & e_3 & e_4 \\
                e_3 & e_4-\frac{1}{4}g_2e_0 & e_5 \\
                e_4 & e_5 & e_6 
                 \endmatrix \right|+\frac{1}{4}g_3\left| \matrix
                0 & e_0 & e_2 \\
                e_0 & e_2 & e_4 \\
                e_2 & e_4 & e_6 
                 \endmatrix \right|$$
  These elements are central in the Poisson algebra $S^*\Cal F_6$.
  
  Let $(f_{\alpha,\beta}),1\leqslant\alpha,\beta\leqslant p+1$ be a $(p+1)\times(p+1)$-matrix of the elements in $\Cal F$.
  Then $\det(f_{\alpha,\beta})\in S^{p+1}\Cal F$ defines some element of degree $p+1$. It follows from our definition of 
  $x_p$, that if rank of $(f_{\alpha,\beta})$ as matrix of functions is equal to 1 (that is 
  $f_{\alpha,\beta}(z)f_{\alpha^{\prime},\beta^{\prime}}(z)=f_{\alpha,\beta^{\prime}}(z)f_{\alpha^{\prime},\beta}(z)$ for 
  usual product of functions), then $\det(f_{\alpha,\beta})\in\ker x_p$ as element in $S^{p+1}\Cal F$. Moreover, let 
  us extend the definition of the functions $e_{\alpha}$ to all integral $\alpha$. Let $\widetilde{\Cal F}$ is the space 
  of functions spanned by $\{e_{\alpha};\alpha\in\Bbb Z\}$. Let $\widetilde{x}_p$ be the natural extension of $x_p$ to
  $S^*\widetilde{\Cal F}$. If $(f_{\alpha,\beta}^{(i)})$ are $(p+1)\times(p+1)$-matrices of the elements from
  $\widetilde{\Cal F}$ such that rank $(f_{\alpha,\beta}^{(i)})$ is equal to 1 for each $i$, then 
  $\det(f_{\alpha,\beta}^{(i)})\in
  \ker\widetilde x_p$ as element from $S^{p+1}\widetilde{\Cal F}$. But if some linear combination 
  $\Psi=\sum_i\lambda_i\det(f_{\alpha,\beta}^{(i)})$ belongs to $S^{p+1}\Cal F$, then $\Psi\in\ker x_p$. For example, 
  the elements $C^{(2)}_1,C^{(3)}_1$ may be written as follow: 
  $$C^{(2)}_1=\left| \matrix
                e_2 & e_3\\
                e_3 &  e_4-\frac{1}{4}g_2e_0-\frac{1}{4}g_3e_{-2}
                 \endmatrix \right|+\frac{1}{4}g_3\left| \matrix
                e_{-2} & e_0\\
                e_0 & e_2
                 \endmatrix \right|,$$
       $$C^{(3)}_1=\left| \matrix
                e_2 & e_3 & e_4 \\
                e_3 & e_4-\frac{1}{4}g_2e_0-\frac{1}{4}g_3e_{-2} & e_5 \\
                e_4 & e_5 & e_6 
                 \endmatrix \right|+\frac{1}{4}g_3\left| \matrix
                e_{-2} & e_0 & e_2 \\
                e_0 & e_2 & e_4 \\
                e_2 & e_4 & e_6 
                 \endmatrix \right|.   $$

        In general case of even 
  $n$ the construction of central elements $C^{(\frac{n}{2})}_0$ and $C^{(\frac{n}{2})}_1$ is similar. Define matrix 
  $(g_{\alpha,\beta}),1\leqslant\alpha,\beta\leqslant\frac{n}{2}$ as follows: $g_{1,1}=e_0, g_{1,\alpha}=g_{\alpha,1}=
  e_{\alpha}$ for $\alpha>1$; $g_{\alpha,\beta}(z)=e_{\alpha}(z)e_{\beta}(z)$ for $\alpha,\beta>1$. For example, 
  $g_{2,2}(z)=\wp(z)^2$, so $g_{2,2}=e_4$ and $g_{3,3}(z)=\frac{1}{4}\wp^{\prime}(z)^2=\wp(z)^3-\frac{1}{4}g_2\wp(z)-
  \frac{1}{4}g_3$, so $g_{3,3}=e_6-\frac{1}{4}g_2e_2-\frac{1}{4}g_3e_0$. We set $C^{(\frac{n}{2})}_0=
  \det(g_{\alpha,\beta})$ as element in $S^{\frac{n}{2}}\Cal F_n$. Define the matrix
  $(g_{\alpha,\beta}^{(1)}),1\leqslant\alpha,\beta\leqslant\frac{n}{2}$ as follow: $g_{1,\alpha}^{(1)}=g_{\alpha,1}^{(1)}=
  e_{\alpha+1}$ and $g_{\alpha,\beta}^{(1)}(z)=\frac{e_{\alpha+1}(z)e_{\beta+1}(z)}{e_2(z)}\in\widetilde{\Cal F}$. It is 
  clear that all $g_{\alpha,\beta}^{(1)}$ belong to $\Cal F$ except $g_{2,2}^{(1)}$. We have $g_{2,2}^{(1)}=e_4-\frac{1}{4}
  g_2e_0-\frac{1}{4}g_3e_{-2}$. Let us define the matrix $(g_{\alpha,\beta}^{(2)}),1\leqslant\alpha,\beta\leqslant\frac{n}{2}$
   as follows: $g_{1,1}^{(2)}=e_{-2}, g_{1,2}^{(2)}=g_{2,1}^{(2)}=e_0, g_{1,\alpha}^{(2)}=g_{\alpha,1}^{(2)}=e_{\alpha-1}$
  for $3\leqslant\alpha\leqslant\frac{n}{2}$; $g_{\alpha,\beta}^{(2)}(z)=\frac{g_{\alpha,1}^{(2)}(z)g_{1,\beta}^{(2)}(z)}
  {g_{1,1}^{(2)}(z)}$. We set $C^{(\frac{n}{2})}_1=\det(g_{\alpha,\beta}^{(1)})+\frac{1}{4}g_3
  \det(g_{\alpha,\beta}^{(2)})$. In fact, $C^{(\frac{n}{2})}_1\in S^{\frac{n}{2}}\Cal F_n$ and $x_{\frac{n}{2}-1}
  (C^{(\frac{n}{2})}_1)=0$.
  
  Let us construct the elements $C^{(n)}\in S^n\Cal F_n$ for odd $n$ such that $C^{(n)}\in\ker x_{\frac{n-1}{2}}$. It is 
  clear that our elements $C^{(\frac{n+1}{2})}_0$ and $C^{(\frac{n+1}{2})}_1$ from $S^{\frac{n+1}{2}}\Cal F_{n+1}$ have 
  a form: $C^{(\frac{n+1}{2})}_0=A_0+B_0e_{n+1}$ and $C^{(\frac{n+1}{2})}_1=A_1+B_1e_{n+1}$ where $A_0,A_1\in 
  S^{\frac{n+1}{2}}\Cal F_n$ and $B_0,B_1\in S^{\frac{n-1}{2}}\Cal F_n$. We set $C^{(n)}=B_0C^{(\frac{n+1}{2})}_1-B_1
  C^{(\frac{n+1}{2})}_0=B_0A_1-B_1A_0$. It is clear that $C^{(n)}\in\ker x_{\frac{n-1}{2}}$. 
  
  {\bf Proposition 4.} {\it The center of the Poisson algebra $S^*\Cal F_n$ is generated by $C^{(\frac{n}{2})}_0$ and
  $C^{(\frac{n}{2})}_1$ for even $n$. For odd $n$ the center is generated by $C^{(n)}$.}
  
  {\bf Proof.} One can check that the elements $C^{(\frac{n}{2})}_0,C^{(\frac{n}{2})}_1$ (resp. $C^{(n)}$) for even 
  (resp. odd) $n$ are central in $S^*\Cal F_n$. On the other hand, the quotient algebra $S^*\Cal F_n\diagup
  (C^{(\frac{n}{2})}_0,C^{(\frac{n}{2})}_1)$ (resp. $S^*\Cal F_n\diagup(C^{(n)})$) is isomorphic to the image of the 
  homomorphism $x_{\frac{n}{2}-1}$ (resp. $x_{\frac{n-1}{2}}$) which is the algebra of functions on the symplectic manifold
  and has a trivial center. So the center is generated by $C^{(\frac{n}{2})}_0,C^{(\frac{n}{2})}_1$ (resp. $C^{(n)}$).
  
  {\bf Remark.} Considering any two linear combinations of our three Poisson brackets one obtains a bihamiltonian structure.
  Lenard scheme provides a family of commuting elements from Casimir elements $C^{(\frac{n}{2})}_0,C^{(\frac{n}{2})}_1$
  for even $n$ and $C^{(n)}$ for odd $n$.

  \newpage
  \centerline{\bf Appendix}
  \centerline{\bf Elliptic functions}
  \medskip
  
  For an integral lattice $\Gamma\subset\Bbb C$ the Weierstrass elliptic function is defined as follows:
  $$\wp(z)=\frac{1}{z^2}+\botshave{\sum_{\omega\in\Gamma^{\prime}}}(\frac{1}{(z-\omega)^2}-\frac{1}{\omega^2}),
  \text{ where }\Gamma^{\prime}=\Gamma\diagdown\{0\}$$
  One has: $\wp^{\prime}(z)^2=4\wp(z)^3-g_2\wp(z)-g_3$, where $g_2$ and $g_3$ depend on the lattice $\Gamma$ only.
  
  The Weierstrass zeta function is defined as follows:
  $$\zeta(z)=\frac{1}{z}+\botshave{\sum_{\omega\in\Gamma^{\prime}}}(\frac{1}{z-\omega}+\frac{1}{\omega}+
  \frac{z}{\omega^2})$$
  The function $\zeta(z)$ is not elliptic, but one has: $\zeta(z+\omega)=\zeta(z)+\eta(\omega)$, where $\eta:\Gamma\to
  \Bbb C$ is a $\Bbb Z$-linear function. The function $\zeta(z_1-z_2)-\zeta(z_1)+\zeta(z_2)$ is elliptic in variables
  $z_1$ and $z_2$. It is clear that $\zeta(-z)=-\zeta(z)$, $\zeta^{\prime}(z)=-\wp(z)$. One has the following useful
   decomposition:
  $$\wp(z)=\frac{1}{z^2}+\frac{1}{20}g_2z^2+\frac{1}{28}g_3z^4+...$$

  We need the following identities:
  $$(\zeta(x-y)-\zeta(x)+\zeta(y))(\wp(x)-\wp(y))=\frac{1}{2}(\wp^{\prime}(x)+\wp^{\prime}(y))\eqno(5)$$
  $$(\zeta(x-y)-\zeta(x)+\zeta(y))(\wp^{\prime}(x)-\wp^{\prime}(y))=2\wp(x)^2+2\wp(x)\wp(y)+2\wp(y)^2-\frac{1}{2}g_2$$
  Proof of these identities is standard: to calculate the decomposition in the neighbourhood of the point $x=y=0$.
  \medskip
  {\bf Acknowledgments.} This work was supported by the grants INTAS-OPEN-00-00055, RFBR-02-01-01015, RFBR-00-15-96579
  
  \newpage
  \centerline{\bf References}
  \medskip
  
  1. Peter D.Lax, Almost periodic solutions of the KdV equation, SIAM Rev. 18(1976), no.3, 351-375.
  
  2. Franco Magri, A simple model of the integrable Hamiltonian equation, Journal of Mathematical Physics 19(1978),
   no.5, 1156-1162.
   
  3. I.M.Gelfand and I.Ja.Dorfman, Hamiltonian operators and algebraic structures associated with them, 
  Funktsional. Anal. i Prilozhen. 13(1979), no.4.
  
  4. I.M.Gelfand and I.Zakharevich, Webs, Lenard schemes and the local geometry of bihamiltonian Toda and Lax structures, 
  math.DG/9903080.
  
  5. B.L.Feigin, A.V.Odesskii, Vector bundles on elliptic curve and Sklyanin algebras. Amer.Math.Soc.Transl.Ser.2.V.185.
  
  6. A.V.Odesskii, Elliptic algebras, Russian Mathematical Surveys, Vol.57(2002),no.6.
  
  7. S.Lang, Elliptic functions, Addson-Wesley Publishing Company, INC., 1973.
  
  \bigskip
  
  LANDAU INSTITUTE FOR THEORETICAL PHYSICS
  
  INSTITUT DES HAUTES ETUDES SCIENTIFIQUES
  
  E-mail address: odesskii\@mccme.ru

 \enddocument